\documentclass[graybox]{svmult}

\usepackage{type1cm}        
\usepackage{makeidx}         
\usepackage{graphicx}        
\usepackage{multicol}        
\usepackage[bottom]{footmisc}
\usepackage{soul}

\newcommand{\bC}[0]{\mathbf{C}}

\usepackage{newtxtext}       %
\usepackage[varvw]{newtxmath}       

\makeindex             

\usepackage{cleveref}
\usepackage{multirow}

\crefname{equation}{}{}  

\begin{document}

\title*{Preconditioning a Fluid--Structure Interaction Problem Using Monolithic and Block Domain Decomposition Methods for the Fluid}
\titlerunning{Preconditioning an FSI Problem Using Monolithic and Block DDM for the Fluid} 
\author{%
Axel Klawonn\orcidID{0000-0003-4765-7387}, 
Jascha Knepper\orcidID{0000-0002-8769-2235}, and 
Lea Saßmannshausen\orcidID{0009-0005-9538-5461}}
\institute{Axel Klawonn, Jascha Knepper, and Lea Saßmannshausen \at Department of Mathematics and Computer Science, and Center for Data and Simulation Science, University of Cologne, \email{[axel.klawonn, jascha.knepper, l.sassmannshausen]@uni-koeln.de }
}
%
%
\maketitle

\abstract{
A fluid–structure interaction (FSI) problem is solved via a monolithic coupling of the fluid, structure, and geometry subproblems. 
The iterative GMRES solver is accelerated with the FaCSI block preconditioner.  
In the FaCSI factorization, the fluid subproblem is approximated using either a monolithic preconditioner or the block preconditioner SIMPLE. 
Two-level overlapping Schwarz methods are then used to approximate the arising inverses. 
The robustness and scalability of the monolithic and SIMPLE preconditioners are compared for a realistic patient-specific artery.
The results indicate that the monolithic preconditioning of the fluid subproblem performs better than the SIMPLE approach. 
Different flow rates are tested and parallel strong scaling has been evaluated.
}


\section{Introduction}
To solve a fluid--structure interaction (FSI) problem, monolithic and segregated approaches are both commonly used, where segregation means that fluid and structure problems are solved separately; see \cite{FSI_overview_2024} for a summary of the different methods and the references therein. 
For the monolithic approach, the FaCSI block preconditioner \cite{FaCSI_DEPARIS2016700} can be used. 
The fluid component in the FaCSI block preconditioner is preconditioned by the SIMPLE block preconditioner. 
In \cite{heinlein2025monolithicblockoverlappingschwarz}, monolithic and block overlapping preconditioners were compared for fluid flow problems modeled by the incompressible Navier--Stokes equations. 
The same geometry of a realistic artery as in \cite{heinlein2025monolithicblockoverlappingschwarz} is used. 
The inverses arising in the FaCSI factorization are approximated by two-level overlapping Scharz preconditioners. 
For comparison, we will replace SIMPLE with a monolithic preconditioner.
The primary focus of this paper is to analyze whether the findings in \cite{heinlein2025monolithicblockoverlappingschwarz} can be applied to the fluid subproblem in the context of an FSI problem.

We will briefly describe the FSI problem, the FaCSI factorization that is used for the monolithic solver, and the preconditioners. 
Let us note that, here, the denominations ``monolithic'' and ``block'' preconditioners refer only to the fluid subproblem and not to the complete coupled FSI problem.
We can observe in \Cref{fig: flowrate and scaling} that the monolithic approach for the fluid subproblem in contrast to the block preconditioning one offers a lower average iteration count and more robustness for high flow rates.
Additionally, we can observe good strong scaling properties; see also \Cref{Table: Scalting details}.


\section{Fluid--Structure Interaction Problem}
\label{sec:1}

For a thorough description of FSI problems as discussed in this section, we refer to  \cite{CardiovascularMathematics,balzani:2015:nmf}. 
Let $\Omega^{t}\subset\mathbb{R}^3$ be the union of the fluid domain $\Omega_{\!f}^{t} \subset\Omega^t$ and the solid domain $\Omega_{\!s}^{t} \subset\Omega^t$ at time $t$. 
The initial or reference configuration is at time $t=0$, that is, $\Omega^{0} = \Omega_{\!f}^{0} \cup \Omega_s^{0}$.
The interface between the fluid and solid domains is defined as $\Gamma = \partial \Omega_{\!f}^{0} \cap \partial \Omega_s^{0}$. 
The solid problem is modeled by an isotropic nonlinear neo-Hookean material model.
The displacement $\mathbf{d}_s$ of the solid domain in the reference configuration is governed by
\begin{align*}
\rho_s \frac{\partial^{2} \mathbf{d}_s}{\partial t ^{2}} - \nabla \cdot (\mathbf{FS}) = \rho_s \mathbf{f}_s ~~\text{in}~\Omega_s^{0} \times (0,T],
\end{align*}
where $\rho_s$ is the volumetric mass density of the solid, and $\mathbf{FS}$ are the first Piola--Kirchhoff stresses.
To model the fluid in the moving domain $\Omega_{\!f}^{t}$, the ALE (arbitrary Lagrangian--Eulerian) map is used.
We define an ALE mapping $\mathcal{A}_f$ from $\Omega_{\!f}^{0}$ to $\Omega_{\!f}^{t}$ via
\begin{align*}
\mathcal{A}_{\!f}: \Omega_{\!f}^{0} \rightarrow \Omega_{\!f}^{t},~~ \mathcal{A}_{\!f}(\mathbf{X}) = \mathbf{X} + \mathbf{d}_{\!f}(\mathbf{X})\quad\forall \, \mathbf{X} \subseteq \Omega_{\!f}^{0},
\end{align*}
where $\mathbf{X}$ are the ALE coordinates, and $\mathbf{d}_{\!f}$ is the displacement of the fluid domain. 
The transient incompressible Navier--Stokes equations -- to model a Newtonian fluid on the moving domain $\Omega_{\!f}^{t}$ -- are then written in the ALE formulation as
\begin{align*}
\rho_{\!f} 
	\Big ( 
	\frac{\partial \mathbf{u}_{\!f}}{\partial t} \Big|_{\mathbf{X}} 
	+  \big[(\mathbf{u}_{\!f} - \mathbf{w}) \cdot \nabla \big] \mathbf{u}_{\!f} 
	\Big ) 
	- \nabla \cdot \sigma_{\!f}(\mathbf{u}_{\!f},p)  
&= 0 
	\quad\text{ in }\quad
	\Omega_{\!f}^{t} \times (0,T], \\
\operatorname{div}(\mathbf{u}_{\!f} ) 
&= 0
	\quad\text{ in }\quad
	\Omega_{\!f}^{t} \times (0,T],
\end{align*}
with the fluid density $\rho_{\!f}$, fluid velocity $\mathbf{u}_{\!f}$, and pressure $p$. 
Then, $\mathbf{w} = \frac{\partial \mathbf{d}_{\!f}}{\partial t}\big|_{\mathbf{X}}$ describes the velocity of the moving fluid domain. 
Note that the time derivative is taken with respect to the ALE coordinates, as indicated by $\mathbf{X}$.
The Cauchy stress tensor is given by $\sigma_{\!f}(\mathbf{u}_{\!f},p)=\mu_{\!f} \, (\nabla \mathbf{u}_{\!f} + (\nabla \mathbf{u}_f)^{T}) - p \mathbf{I}$, where $\mu_f$ is the dynamic viscosity.

The displacement $\mathbf{d}_{\!f}$ of each point in the fluid domain with respect to the reference configuration $\Omega_{\!f}^{0}$ is determined by a harmonic extension of the moving fluid--structure interface, i.e., by solving the problem
\begin{equation}
\label{eq:geometryProblem}
\Delta \mathbf{d}_{\!f} = 0 \text{ in } \Omega_{\!f}^{0},
\quad 
\text{with }\mathbf{d}_{\!f} = \mathbf{d}_s \text{ on }\Gamma,
\quad
\text{and } \mathbf{d}_{\!f} \cdot \mathbf{n}_{\!f} = 0 \text{ on }\partial \Omega_{\!f}^{0} \setminus \Gamma.
\end{equation}
$\mathbf{n}_{\!f}$ is the outward normal vector field of the fluid domain; 
$\mathbf{n}_s$ is the corresponding field for the solid domain. 
We refer to \cref{eq:geometryProblem} as the geometry subproblem. 
The geometric adherence of the fluid domain and solid domain is given by $\mathbf{d}_{\!f} = \mathbf{d}_s ~\text{on}~\Gamma \times (0,T]$, and the coupling conditions for the velocities and the stresses on the interface are
\begin{align*}
\frac{\partial \mathbf{d}_s}{\partial t} = \mathbf{u} \circ \mathcal{A}_{\!f} , ~\text{and}~~~ (\text{det}[\mathbf{F}])^{-1} \mathbf{F}^{-T}\! \sigma_{\!f} \mathbf{n}_{\!f} \circ \mathcal{A}_{\!f} + (\mathbf{F} \mathbf{S}) \mathbf{n}_s = 0.
\end{align*}
The finite element method is used for the spatial discretization of the FSI problem. 
The domain is partitioned into tetrahedral elements. 
Here, we consider a P2--P1 discretization for the fluid subproblem and a P2 discretization for the solid subproblem.


\section{The FaCSI Block Preconditioner for an FSI Problem}
The name of the block preconditioner originates from \textbf{Fa}ctorized FSI block system, with static \textbf{C}ondensation and use of the \textbf{SI}MPLE preconditioner for the fluid problem and was introduce in \cite{FaCSI_DEPARIS2016700} by Deparis, Forti, Grandperrin, and Quarteroni.
The FaCSI block preconditioner can be combined with different strategies for the approximation of the fluid subproblem; see  \cite{diss_hochmuth_2020} for a combination with a monolithic two-level overlapping Schwarz preconditioner. 
We will give a brief definition of FaCSI. 
In \cite{FaCSI_DEPARIS2016700} a detailed derivation can be found.
The preconditioner is based on an incomplete block factorization. 
Starting from the Jacobian matrix of the FSI system, which consists of the discretized fluid subproblem $\mathcal{F}$, solid subproblem $\mathcal{S}$, geometry problem $\mathcal{G}$, coupling blocks $\bC_i$, and shape derivatives $\mathcal{D}$, the block preconditioner 
\begin{align}
\label{eq: FaCSI decomp}
\begin{pmatrix}
\mathcal{S} & 0 & 0 & \bC_4 \\
\bC_5 & \mathcal{G} & 0 & 0 \\
0 & \mathcal{D} & \mathcal{F} & \bC_3 \\
\bC_2 & 0 & \bC_1 & 0 
\end{pmatrix}
\approx 
\underbrace{
\begin{pmatrix}
\mathcal{S} & 0 & 0 &0 \\
0 & I_{\mathcal{G}} & 0 & 0 \\
0 & 0 & I_{\mathcal{F}} & 0 \\
0 & 0 &0 & I_\Gamma 
\end{pmatrix}
}
_{\mathcal{B}_{\mathcal{S}}}
\underbrace{
\begin{pmatrix}
I_{\mathcal{S}} & 0 & 0 &0 \\
\bC_5 &\mathcal{G} & 0 & 0 \\
0 & 0 & I_{\mathcal{F}} & 0 \\
0 & 0 &0 & I_\Gamma 
\end{pmatrix}
}
_{\mathcal{B}_{\mathcal{G}}}
\underbrace{
\begin{pmatrix}
I_{\mathcal{S}} & 0 & 0 &0 \\
0 & I_{\mathcal{G}} & 0 & 0 \\
0 & \mathcal{D} & \mathcal{F} & \bC_3 \\
\bC_2 & 0 & \bC_1 & 0
\end{pmatrix}
}
_{\mathcal{B}_{\mathcal{F}}}
=
\mathcal{B}_{\text{FaCSI}}
\end{align}
is defined. 
The approximation in \cref{eq: FaCSI decomp} is obtained by decoupling the different physical problems and neglecting the block $\bC_4$. 
Due to its block-diagonal structure, for an approximation of the inverse of $\mathcal{B}_{\mathcal{S}}$, only an approximation of $\mathcal{S}^{-1}$ is required. 
The block components $\mathcal{B}_{\mathcal{F}}$ and $\mathcal{B}_{\mathcal{G}}$ are factorized further, until the inverses can be approximated efficiently. 
$\mathcal{B}_{\mathcal{G}}$ can easily be split in two factors; see \cite[eq.~(19)]{FaCSI_DEPARIS2016700}.
The fluid component $\mathcal{B}_{\mathcal{F}}$ is factorized as
\begin{align*}
\mathcal{B}_{\mathcal{F}}=
\begin{pmatrix}
I_{\mathcal{S}} & 0 & 0 &0 \\
0 & I_{\mathcal{G}} & 0 & 0 \\
0 & \mathcal{D} & I_{\mathcal{F}} & 0 \\
0 & 0 & 0 & I_\Gamma
\end{pmatrix}
\begin{pmatrix}
I_{\mathcal{S}} & 0 & 0 &0 \\
0 & I_{\mathcal{G}} & 0 & 0 \\
0 & 0 & I_{\mathcal{F}} & 0 \\
\bC_2 & 0 & 0 & I_\Gamma
\end{pmatrix}
\begin{pmatrix}
I_{\mathcal{S}} & 0 & 0 &0 \\
0 & I_{\mathcal{G}} & 0 & 0 \\
0 & 0 & \mathcal{F} & \bC_3 \\
0 & 0 & \bC_1 & 0
\end{pmatrix}
=
\mathcal{B}_{\mathcal{F}_1}\mathcal{B}_{\mathcal{F}_2}\mathcal{B}_{\mathcal{F}_3}.
\end{align*} 
The factors $\mathcal{B}_{\mathcal{F}_1}$ and $\mathcal{B}_{\mathcal{F}_2}$ can be inverted exactly.
The following type of $2\times 2$ block system in $\mathcal{B}_{\mathcal{F}_3}$ remains to be solved:
\begin{align}\label{eq:FaCSI fluid subproblem}
\begin{pmatrix}
\mathcal{F} & \bC_3 \\
\bC_1 & 0 
\end{pmatrix}
\begin{pmatrix}
x \\
\lambda
\end{pmatrix}
= 
\begin{pmatrix}
r_{\mathcal{F}} \\
r_{\lambda}
\end{pmatrix}.
\end{align}
Since $C_1$ restricts to the velocity interface degrees of freedom, and since $C_3=C_1^T$, one can show that $\lambda$ can be eliminated, and a condensed system of the type 
\begin{align}\label{eq: FaCSI fluid problem}
\mathcal{F}_{II} x_{I} = r_{I} - \mathcal{F}_{I\Gamma} x_{\Gamma}
\end{align}
needs to be solved, where $I$ and $\Gamma$ denote the interior and interface degrees of freedom, respectively.
The inverse of $\mathcal{F}_{II}$ in \Cref{eq: FaCSI fluid problem} can now be approximated with a monolithic or block preconditioner.
For the approximation of the other arising inverses, two-level overlapping Schwarz preconditioners are used; see \Cref{sec:schwarz}.


\subsection{The SIMPLE Block Preconditioner for a Fluid Problem}
The SIMPLE (Semi-Implicit Method for Pressure-Linked Equations) preconditioner was originally introduced by Patankar and Spalding \cite{Patankar_SIMPLE_1972} to solve the Navier–Stokes equations; 
the SIMPLE block preconditioner is based on this algorithm.
Derived from a block factorization of the saddle-point problem, the preconditioner is constructed as follows
\begin{equation*}
\mathcal{B}_{\text{SIMPLE}} = 
\begin{bmatrix}
F & 0 \\
B & S_{\text{SIMPLE}} 
\end{bmatrix}
\begin{bmatrix}
I & \frac{1}{\alpha} H_F B^{T} \\
0 & \frac{1}{\alpha} I 
\end{bmatrix},
\end{equation*}
with the Schur complement $S=-C-B F B^T$ of the Navier--Stokes system approximated by $S_{\text{SIMPLE}}=-C-B H_F B^{T}$. 
The diagonal matrix $H_F$ depends on the specific SIMPLE variant. 
We denote the matrix of the default variant SIMPLE as $H_F^D$ and the one of SIMPLEC as $H_F^{\Sigma}$:
\begin{align*}
H_F^D = \operatorname{diag}(F)^{-1} ~ \text{(SIMPLE)}~,~~~
H_F^{\Sigma} = \delta^{ij}\Big(\sum_{k=1}^{N_u} |F_{i,k}|\Big )^{-1}~ \text{(SIMPLEC)}.
\end{align*}


\section{Additive Overlapping Schwarz Preconditioners} \label{sec:schwarz}

We will solve the discretized FSI problem in~\Cref{sec:1} using a Newton--Krylov approach. More precisely, we solve the linear tangent system in each Newton step using a preconditioned GMRES method.
In our implementation, we use the \texttt{Trilinos} \cite{trilinos} package \texttt{FROSch} \cite{Heinlein:2020:FRO} for the Schwarz preconditioners. 
A detailed description of all libraries and packages involved can be found in \cite[Sec. 5]{heinlein2025monolithicblockoverlappingschwarz}.

We give a brief introduction to two-level overlapping Schwarz preconditioners, which are used to precondition the subproblems arising in \cref{eq: FaCSI decomp}. Since the fluid subproblem is a 2$\times$2 saddle-point system, it can either be preconditioned by a monolithic or block preconditioner. 
In both cases, overlapping Schwarz preconditioners are used to approximate the arising inverses; see \cite[Sec. 4]{heinlein2025monolithicblockoverlappingschwarz} for a detailed description.
Let $\Omega$ be an arbitrary domain and let it be decomposed into nonoverlapping subdomains $\{ \Omega_i \}_{i=1}^{N}$.  To construct the overlapping subdomains we extend them by $k$ layers of finite elements, resulting in an overlapping domain decomposition $\{ \Omega_i' \}_{i=1}^{N}$. 
Based on the overlapping domain decomposition, we define a restriction operator $R_i: V^h \rightarrow V_i^h$, $i=1,\dots,N$, to map from the global finite element space $ V^h(\Omega)$ to the local finite element space $V_i^h(\Omega)$ on the overlapping subdomains $\Omega_i'$. 
The corresponding prolongation operator $R_i^T$ is the transpose of the restriction operator and extends a local function on $\Omega_i'$ by zero outside of $\Omega_i'$. 
Using both operators, $R_i$, $R_i^{T}$, we can define local overlapping stiffness matrices $K_i := R_i K R_i^{T}$, $i=1,\ldots,N$. 
The one-level Schwarz preconditioner does in general not scale with the number of subdomains, since the number of Krylov iterations increases with a growing number of subdomains. 
To define a numerically scalable overlapping Schwarz preconditioner, which yields a convergence rate independent of the number of subdomains, a coarse level needs to be introduced. 
With a coarse interpolation operator $\Phi:V_0\to V^h$, the additive two-level overlapping Schwarz preconditioner reads
\begin{equation} \label{eq:two_level schwarz}
	M_{\text{OS-2}}^{-1}
	=
	\Phi K_0^{-1} \Phi^{T}
	+
	M^{-1}_{\text{OS-1}},
	\quad
	M^{-1}_{\text{OS-1}} = \sum_{i=1}^{N} R_i^T K_i^{-1} R_i.
\end{equation}
Here, $K_0 = \Phi^{T} K \Phi$ is the coarse matrix, which is a Galerkin projection of~$K$ into the coarse space. 
The columns of $\Phi$ are a basis of the coarse space $V_0$. 
In this paper, we use GDSW-type coarse spaces, which originate from \cite{dohrmann_domain_2008, dohrmann_design_2017, heinlein2025monolithicblockoverlappingschwarz}.
Additional insights and more in-depth explanations can be found in related work:
The first comparison of monolithic and block overlapping Schwarz preconditioners for saddle-point problems was made in \cite{Klawonn_2000_Mono_Block}.
In \cite{heinlein_reduced_2019, heinlein2025monolithicblockoverlappingschwarz}, two-level overlapping Schwarz preconditioners are considered and compared for fluid flow problems. 
In \cite{heinlein_parallel_2016a}, two-level overlapping Schwarz preconditioners are investigated for the solid subproblem in FSI.


\section{Results}
\label{sect:results}
The parallel results in this section were obtained on the Fritz supercomputer at Friedrich-Alexander-Universität Erlangen-Nürnberg. 
The following solver setting with respect to Newton and GMRES were used: 
Newton's method terminates when the relative residual or the Newton update reach a tolerance of $10^{-8}$. 
Newton's method is used with GMRES along with an adaptive forcing term $\eta_k$; see \cite{heinlein2025monolithicblockoverlappingschwarz} for more details. 
The value of $\eta_{k,\min}$ is set to $10^{-4}$, and $\eta_{k,\max}$ is $10^{-8}$.
The initial mesh partition defining the nonoverlapping subdomains is constructed by METIS, which results in an unstructured partition of the mesh.
We use an overlap of $\delta = 1$.

\begin{figure}[t]
\includegraphics[scale=.16]{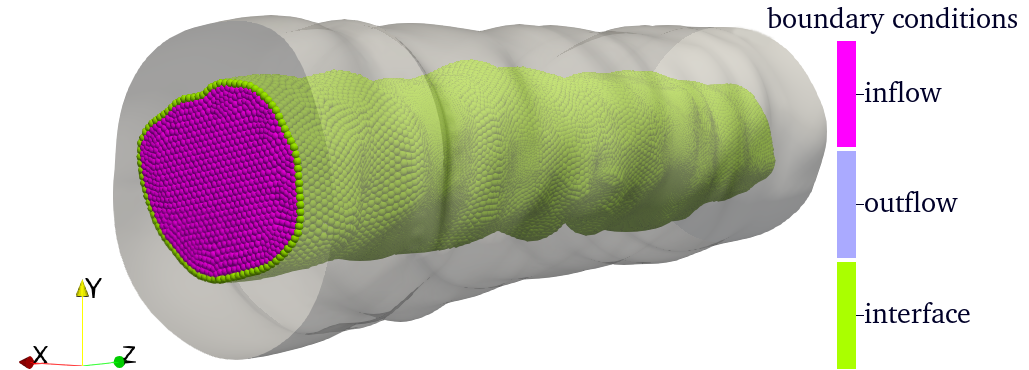}
\includegraphics[scale=.16]{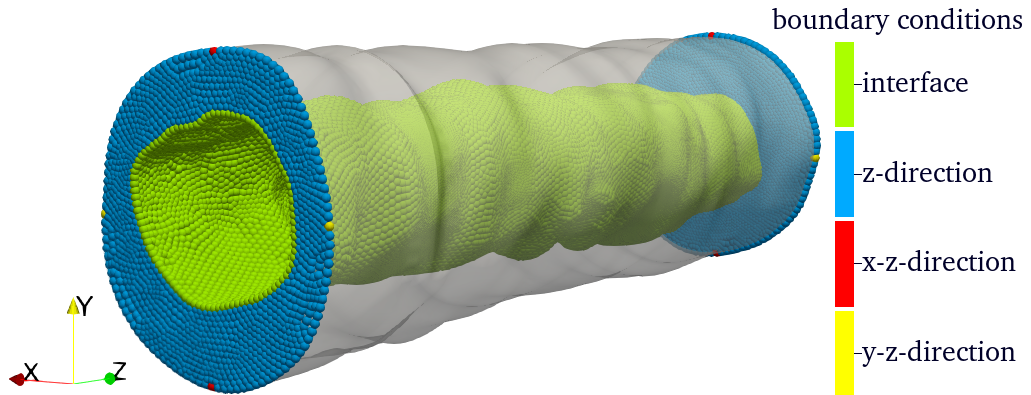}
\caption{Boundary conditions of fluid and solid domain. 
For the fluid subproblem, a parabolic-like inflow profile is prescribed on the inlet (in pink). 
On the outflow, a pressure boundary condition is prescribed.
For the solid subproblem, the $z=0$ and $z=L$ planes are held in $z$ direction. 
Additional points on the respective inflow and outflow planes are held in $x$-$z$ or $y$-$z$ direction to prevent rigid body motions and to ensure that the system is statically determinate.}
\label{fig:1 artery}       
\end{figure}

We consider the realistic, patient specific artery in \Cref{fig:1 artery}, which was obtained for the work in~\cite{Artery_Balzani_2012}. 
The interior (lumen) diameter of the artery varies between 0.2\,cm and 0.34\,cm. 
The solid component contains layers of plaque, media, and adventitia. 
Here, all components are modeled homogeneously.
The outer diameter varies between 0.44\,cm and 0.47\,cm. 
More details can be found in \cite[Sec. 6.1.2.]{heinlein2025monolithicblockoverlappingschwarz}. 
In \Cref{Table: FSI Parameter}, the full list of FSI parameters is shown.

We impose a boundary condition on the outlet of the fluid domain to prescribe a pressure value -- the absorbing boundary condition; 
a detailed description can be found in \cite[Sec. 4.1.5.]{balzani:2015:nmf}. 
We denote the pressure in the steady state, after the initial ramp phase, as the reference pressure $p_{\text{ref}}$.
For the time discretization of the fluid subproblem, we choose a BDF-2 time stepping scheme, and for the solid subproblem a Newmark scheme.
The time step size is $\Delta t= 0.001$\,s.

\begin{figure}[t]
\includegraphics[scale=.55]{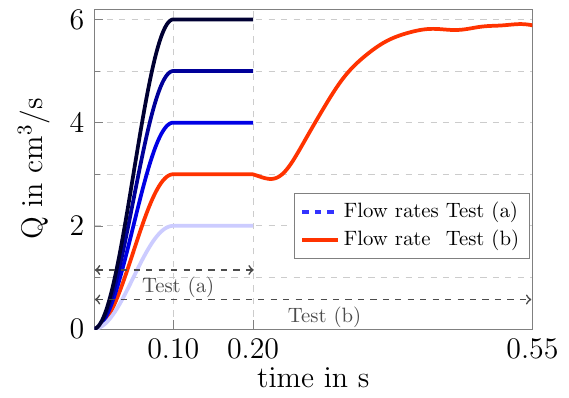}
\includegraphics[scale=.55]{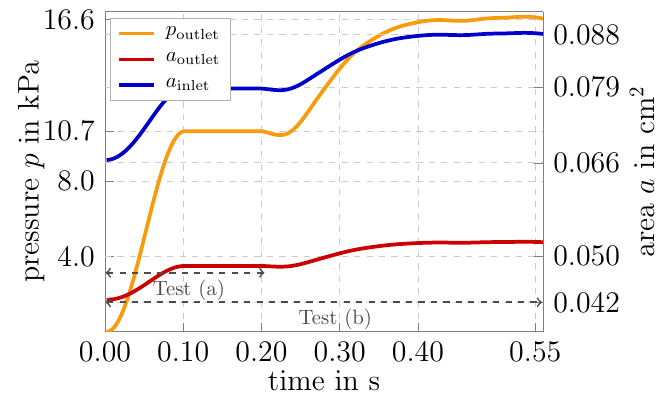}
\caption{\textbf{Left:} Prescribed flowrate on inlet of fluid domain. We consider a range from 2 to 6\,cm$^{3}$/s for Test (a). 
For the steady phase up to 0.2\,s the pressure is not influenced by the flow rate due to the absorbing boundary condition.
\textbf{Right:} The prescribed reference pressure $p_{\text{ref}}$ leads to the pressure measured  at the outlet. It follows the fluid flow rate ramp prescribed (left). 
For Test (b) the flow rate prescribed by the heart beat also influences the pressure value after the steady phase.
}
\label{fig: results pressure and area}
\end{figure}

We consider two different test settings. 
First, we compare results for different flow rates. 
The flow rates vary between 2 and 6\,cm$^{3}$/s. 
The test contains a ramp and a steady phase, each 0.1\,s long; cf. \textit{Test (a)} in \Cref{fig: results pressure and area}.
Secondly, for a strong scaling test, we also include the flow profile of half a heartbeat until 0.56\,s, which is at the peak of the heart beat flow; cf. \textit{Test (b)} in \Cref{fig: results pressure and area} for the resulting pressure, area, and flow rates.
 
\begin{table}[!t]
\centering
\caption{Parameters for the fluid and solid. 
The reference fluid pressure $p$ is kept constant after an initial ramp phase. 
For the flowrate at the end of the initial ramp phase, $Q$, see \Cref{fig: results pressure and area}.} \label{Table: FSI Parameter}
\begin{tabular}{p{1.5cm}p{1.5cm}p{1.7cm}p{1.7cm}p{1.7cm} p{1.7cm}p{1.1cm}}
\hline\noalign{\smallskip}
$\nu$ in cm$^{2}$/s	 &  $ p_{\text{ref}}$ in kPa &  $\rho_f$ in kg/cm$^{3}$  &  $\rho_s$ in kg/cm$^{3}$ & Poisson ratio &  $\mu$ in kPa & $E$ in kPa   \\
\noalign{\smallskip}\svhline\noalign{\smallskip}
0.0291        & 				10.66		&	 1.03\,10$^{-3}$		& 1.0\,10$^{-3}$ & 0.49 	& 		127.52				&	 380			\\
\noalign{\smallskip}\hline\noalign{\smallskip}
\end{tabular}
\end{table}

For preconditioning the fluid subproblem, we will apply findings from \cite{heinlein2025monolithicblockoverlappingschwarz}. 
The coarse space combination GDSW$^{*}$ and RGDSW is used for the velocity and pressure components in the monolithic preconditioner of the fluid problem. 
In case of the SIMPLEC variant of the SIMPLE block preconditioner, the RGDSW coarse space is used for the fluid and Schur complement approximations; cf. \cite[Sec. 4.2.]{heinlein2025monolithicblockoverlappingschwarz}.

\begin{figure}[t]
\hspace*{\fill}
\includegraphics[scale=.55]{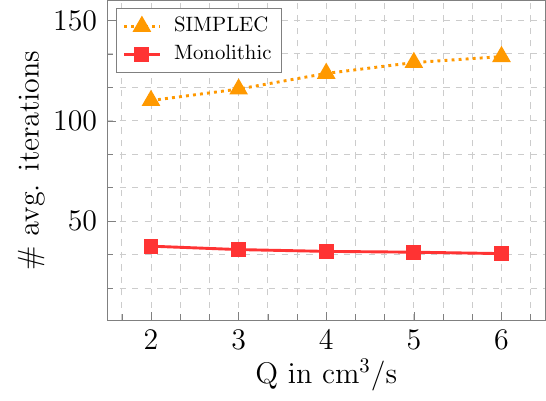}
\hspace*{\fill}
\includegraphics[scale=.55]{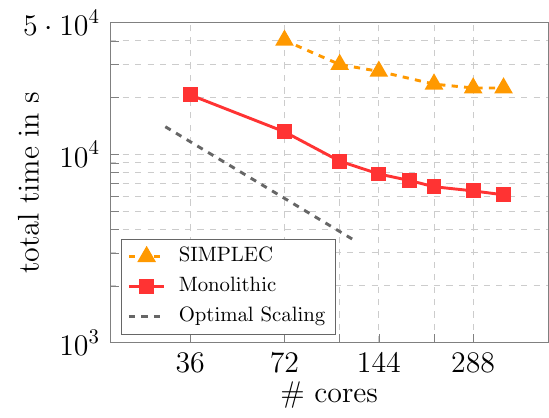}
\hspace*{\fill}
\caption{ Comparison of monlithic and block preconditioner. \textbf{Left:} Resulting average iteration count per Newton step per time step for varying inflow flow rate using Test (a) in \Cref{fig: results pressure and area} for 200 time steps. On average 6 Newton iterations per time step are needed.
\textbf{Right:} Strong scaling using Test (b) in \Cref{fig: results pressure and area} for 600 time steps. The full FSI system has 2\,230\,418 degrees of freedom. Total time consists of setup and solve time. On average 8 Newton iterations per time step are needed.
}
\label{fig: flowrate and scaling}    
\end{figure}

In \Cref{fig: flowrate and scaling}, results for the comparison of block and monolithic preconditioners are shown. 
The figure on the left shows the average iteration count per Newton step per time step to solve the full FSI system. 
The iteration count for the monolithic case is lower compared to the SIMPLEC variant and appears more robust with an increasing flow rate. 
The figure on the right shows the strong scaling behavior of both methods. 
Since the monolithic preconditioner for the fluid component performs better than its block counterpart, the overall iteration count to solve the FSI problem is lower.
In \Cref{Table: Scalting details}, detailed timings for a varying number of processor cores can be found. 
The SIMPLEC block preconditioner offers lower setup times compared to the monolithic approach, but results in significantly longer times for the iterative solvers.
 
\begin{table}[!t]
\caption{Detailed results for the strong scaling plot shown in \Cref{fig: flowrate and scaling} (right). Setup time consists of construction of preconditioner and solve time consists of linear and nonlinear solve. Simulation of 600 time steps with flow rate 3\,cm$^{3}$/s. 
} \label{Table: Scalting details}
\renewcommand{\arraystretch}{1.2}
\begin{tabular}{ |p{1.0cm}| p{1.75cm} p{1.6cm}p{1.6cm} |p{1.75cm}p{1.6cm}p{1.6cm}|}
\hline  
Prec.	&	 \multicolumn{3}{c|}{Monolithic} 		& \multicolumn{3}{c|}{SIMPLEC} 	    	\\ \hline  \hline 
N 	   & 	\# avg. iter. & setup 	& solve 		& \# avg. iter. & setup & solve   \\ 
\hline 
36		&  		24.5	& 	3\,760\,s  &  16\,900\,s	 & 		~~~~-				&    		~~~~-		&   ~~~~- \\				
72		&  		29.3	& 	1\,911\,s	 & 11\,250\,s	 & 	\phantom{0}95.5	&  1\,738\,s 	& 38\,520\,s    \\
144		&   		 30.1	& 	1\,111\,s	 & \phantom{0}6\,746\,s	 & 	108.1		&  \phantom{0}\,965\,s		&  26\,590\,s  \\
288			&   	35.7	 					& 	\phantom{0}883\,s	 &   	\phantom{0}5\,506\,s	 & 	129.2					&    	\phantom{0}\,681\,s			& 21\,800\,s   \\
\hline
\end{tabular}
\end{table}

\vspace{1cm}

{\small
\textbf{Acknowledgments}
Financial funding from the Deutsche Forschungsgemeinschaft (DFG) through the Priority Program 2311 "Robust coupling of continuum-biomechanical in silico models to establish active biological system models for later use in clinical applications - Co-design of modeling, numerics and usability", project ID 465228106, is greatly appreciated.
The authors gratefully acknowledge the scientific support and HPC resources provided by the Erlangen National High Performance Computing Center (NHR@FAU) of the Friedrich-Alexander-Universität Erlangen-Nürnberg (FAU) under the NHR project k105be. NHR funding is provided by federal and Bavarian state authorities. NHR@FAU hardware is partially funded by the German Research Foundation (DFG) -- 440719683.

}
\bibliographystyle{spmpsci.bst}
\bibliography{bib}

\end{document}